\documentclass[hidelinks,onefignum,onetabnum]{siamart220329}

\usepackage{amsfonts,upquote}
\usepackage{graphicx,verbatim}

\headers{$L^2$ and $L^\infty$ rational approximation}{Ackermann, Reiter and Trefethen}

\title{{\boldmath$L^2$ and $L^\infty$ rational approximation on the unit disk}%
\thanks{Submitted to the editors DATE.}}

\author{Michael S. Ackermann\thanks{Dept.\ of Mathematics,
Virginia Tech, Blacksburg, VA 24061, USA
(\email{amike98@vt.edu})}
\and Sean Reiter\thanks{Courant Institute of Mathematical Sciences,
New York University, New York, NY 10012, USA
(\email{s.reiter@nyu.edu})}
\and Lloyd N. Trefethen\thanks{School of Engineering and Applied Sciences,
Harvard University, Cambridge, MA 02138, USA}
(\email{trefethen@seas.harvard.edu})}

\def\complex{{\mathbb{C}}}
\def\Rn{R_n}\def\Rnu{R_\nu}
\def\Rnm{R_{n-1}}
\def\rt{r_2^*}\def\ri{r_\infty^*}
\def\Et{E_2^{}}\def\Ei{E_\infty^{}}
\def\et{e_2^*}\def\ei{e_\infty^*}

\def\rtil{\tilde r}
\def\mk{m_k^{}}\def\zk{\zeta_k^{}}
\def\pk{\pi_k^{}}\def\pz{\pi_1^{}}
\def\pput(#1,#2)#3{\noindent\smash{\raise#2pt\hbox to 0pt
   {\kern #1pt #3\hss}}\ignorespaces}

\def\R{\ensuremath{\mathbb{R}}}
\def\RNN{\R^{N \times N}}
\def\RNp{\R^{N \times p}}
\def\RNq{\R^{N \times q}}
\def\Rqp{\R^{q \times p}}

\def\Rqp{\R^{q \times p}}

\def\rtil{\tilde r}

\def\A{\mathbf A}
\def\B{\mathbf B}
\def\C{\mathbf C}
\def\D{\mathbf D}
\def\E{\mathbf E}

\def\H{\mathbf H}
\def\L{\mathbf L}
\def\M{\mathbf M}
\def\Y{\mathbf Y}

\def\one{\ensuremath{\mathbf 1}}

\begin{document}

\maketitle


\begin{abstract}
Using recently developed algorithms, we compute and compare
best $L^2$ and $L^\infty$ rational approximations of analytic functions on
the unit disk.  Although there is some theory for these problems
going back decades, this may be the first computational study.
To compute the $L^2$ best approximations, we employ
a new formulation of TF-IRKA in barycentric form.
\end{abstract}

\begin{keywords}
rational approximation, IRKA, AAA
\end{keywords}

\begin{MSCcodes}
41A20, 65D15
\end{MSCcodes}

\section{\label{intro}Introduction}

Rational approximation is an old subject, going back to the 19\kern .2pt th
century, but computations can be challenging.  In this paper we
compute and compare best $L^2$ and $L^\infty$ rational approximations
of scalar analytic functions on the unit disk $D$ in the complex plane
$\complex$.  So far as we are aware, no such comparisons have been
published before.  The algorithms we rely on for our computations
are variants of TF-IRKA
for $L^2$ \cite{bg} and AAA-Lawson for $L^\infty$~\cite{lawson}.
They are not fail-safe, but they succeed in many cases.

Let $n\ge 0$ be an integer and let $\Rn$ be the set
of rational functions of degree $n$ that are analytic in the closed
unit disk $\overline{D}$, meaning functions that can be written
in the form $r(z) = p(z)/q(z)$ where $p$ and $q$ are polynomials of degree
(at most) $n$ and all the roots $\{\pk\}$ of $q$ satisfy $|\pk|>1$.
Our $2$- and $\infty$-norms are defined by
\begin{displaymath}
\|f\|_2^{} = \Bigl(\kern 1pt {1\over 2\pi} \kern -1pt\int_S |f(z)|^2 |dz| \Bigr)^{1/2}, \quad
\|f\|_\infty^{} = \sup_{z\in S} |f(z)|, 
\end{displaymath}
where $S$ is the unit circle.  Note that the division by
$2\pi$ makes $\|\cdot\|_2^{}$
the root-mean-square norm, implying $\|f\|_2^{} \le \|f\|_\infty^{}$ for any $f$. 
We will be concerned with approximation of functions $f$ analytic in $\overline{D}$,
so by the maximum modulus principle, approximation on the circle 
implies approximation in the disk. Our functions belong
to the spaces $L^2(S)$, $L^\infty(S)$, $H^2(D)$, and $H^\infty(D)$ as usually defined,
and the problems we address would often be described as instances of $H^2$ or
$H^\infty$ approximation.
Given a function $f$ analytic in the closed unit disk,
we let $\rt$ and $\ri$ denote degree $n$ {\em best approximations\/} defined by
\begin{displaymath}
\|f-\rt\|_2^{} = \inf_{r\in\Rn} \|f-r\|_2^{}, \quad
\|f-\ri\|_\infty^{} = \inf_{r\in\Rn} \|f-r\|_\infty^{}.
\end{displaymath}
It is known that these approximations exist but need not be unique;
see section 3.
(Best approximations in $L^\infty$ are also known as Chebyshev approximations.)
When we write $\rt$ or $\ri$, it is to be understood that this refers
to any best approximation if there are several.
We denote the optimal errors by 
\begin{displaymath}
\Et = \|f-\rt\|_2^{}, \quad
\Ei = \|f-\ri\|_\infty^{}.
\end{displaymath}
Note that $\rt$, $\ri$, $\Et$, and $\Ei$ all depend on $f$ and $n$,
although our notation does not indicate this.

Instead of the closed unit disk, one could equivalently
consider approximation on the closed {\em exterior\/}
of the disk by assuming $f$ is analytic there, including
at $\infty$.  This is the usual formulation in discrete-time
systems theory \cite{Antbook,abgbook,bunse,grimm}, and it is
in some respects more elegant mathematically, including for
investigating the distribution of poles of approximants, which
accumulate as $n\to\infty$ on curves defined by orthogonality
or energy-minimization conditions~\cite{bsy}.  By a M\"obius
transformation, one could go further and transplant $f$ and $r$
to any other disk or disk-exterior or half-plane, but whereas
$L^\infty$ approximation is invariant under such transplantations,
$L^2$ approximation is not.  In particular, the problem of $L^2$
approximation in a half-plane ($H^2$ approximation) requires
the interpolation condition $r(\infty) = f(\infty)$ in order for
the approximation error to be finite.  Half-plane domains are the
natural setting for the model order reduction and reduced order modeling
of continuous-time dynamical
systems, which is the context in which our $L^2$ algorithms have
most often been applied; see section~5.  We mention that recent
work generalizes some of these ideas to regions other than disks,
disk exteriors, and half-planes~\cite{borghi24,borghi}.

We begin in sections 2--4 with examples of best $L^2$ and $L^\infty$
approximations and discussion of associated mathematical properties.
Sections 5 and 6 then discuss algorithms for $L^2$ approximation,
first based on the state-space formulation of the rational
approximant that is so familiar in
the model order reduction and reduced order modeling community,
then based on a barycentric form of the Hermite rational approximant,
for which there is less previous literature.
We will not give details of our $L^\infty$ algorithm, because it
consists of applying the AAA-Lawson method exactly as described
in \cite{lawson}.  For the examples of this paper we sampled at
200 roots of unity, a number which could of course be increased
for more complicated problems.  Alternatively, one could equally
well use adaptive sampling by the continuum AAA-Lawson algorithm
of \cite{continuum}.  The AAA-Lawson algorithm, like the AAA
algorithm it is derived from (which computes near-best rather
than best approximations), does not work with poles, but is just
a greedy descent method involving support points and barycentric
coefficients based on samples of $f$ on the unit circle.  Poles of
the rational approximant are only computed after the fact, once
good support points and corresponding optimal coefficients have
been found.  Our $L^2$ algorithms presented in sections~5 and~6,
by contrast, do work with poles, exploiting the optimality conditions given
in Theorem~2 below, and they adaptively sample $f$
and its derivative at points in the disk rather than on the circle.

\section{\label{ex1} Examples}
\begin{figure}
\begin{center}
\includegraphics[trim=30 65 20 10, clip, scale=.93]{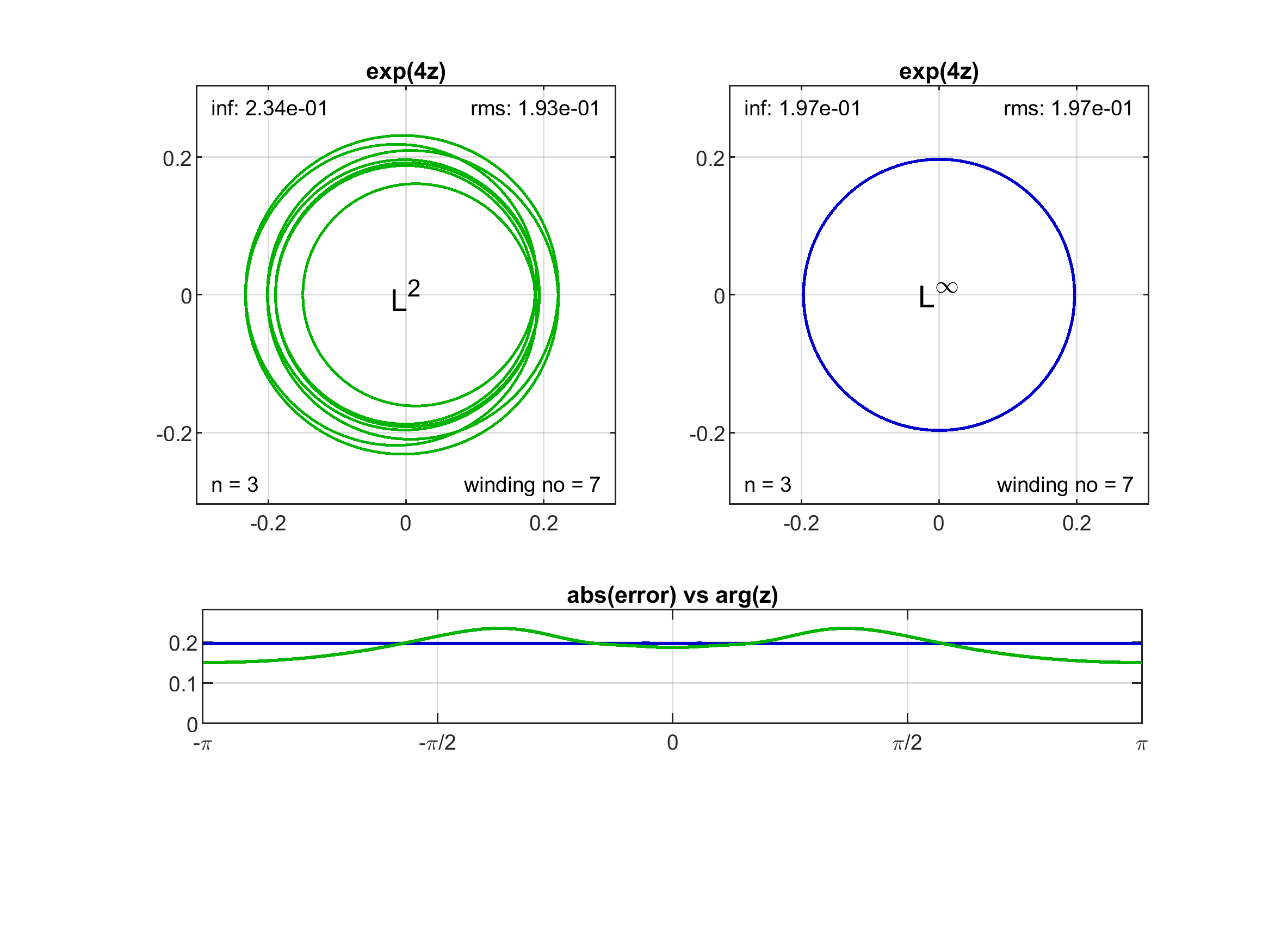}
\end{center}
\caption{\label{fig1}Error curves for
best\/ $L^2$ and\/ $L^\infty$ approximations of degree\/ $n=3$
to $f(z) = \exp(4z)$.  Both curves have winding number\/ $\omega
= 2\kern .3pt n+1=7$.  The
error curve for the\/ $L^\infty$ case is not exactly circular but nearly so.}
\end{figure}

Given $f$ and an approximation $r\in\Rn$, the {\em error function\/} is
\begin{displaymath}
e(z) = r(z) - f(z),
\end{displaymath}
which is defined at least on the closed unit disk, and the {\em
error curve\/} is the closed curve $e(S)$.  The error functions for
the best degree $n$ approximations are denoted by $\et$ and $\ei$.
Figure~\ref{fig1}, the first of several of our figures in the same
format, plots error curves for best degree $3$ approximations of
$f(z)= \exp(4z)$.  Here and in all our examples, except as indicated
otherwise, we believe (usually without proof) that the approximations
we plot are unique best approximations correct to plotting accuracy.
Throughout this paper, we plot curves for $L^2$ best approximations
in green and those for $L^\infty$ best approximations in blue.

\begin{figure}
\begin{center}
\includegraphics[trim=30 65 20 10, clip, scale=.93]{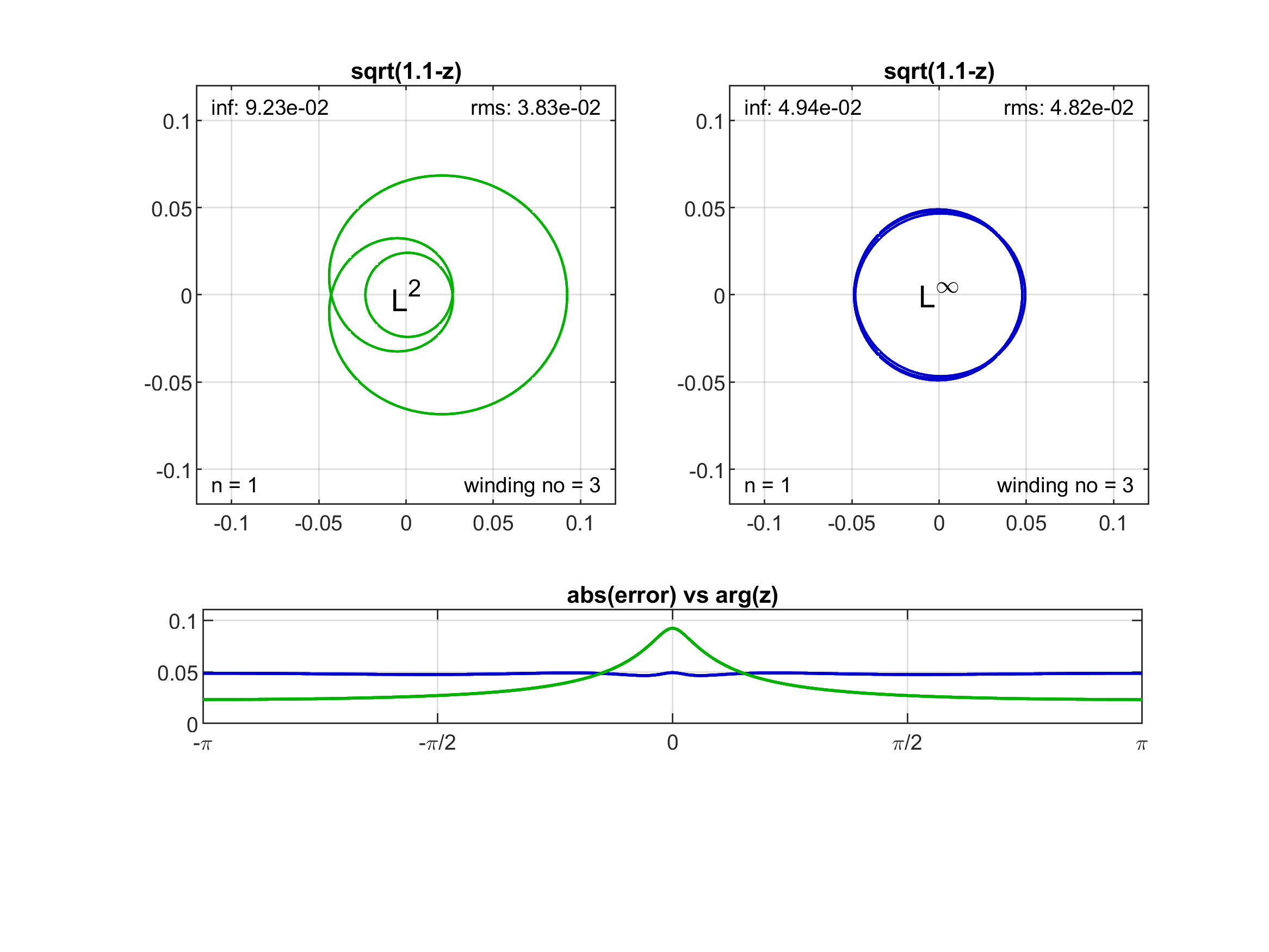}
\end{center}
\caption{\label{fig2}Degree\/ $1$ approximations of\/ $\sqrt{1.1-z}$.}
\end{figure}

\begin{figure}
\begin{center}
\includegraphics[trim=30 65 20 10, clip, scale=.93]{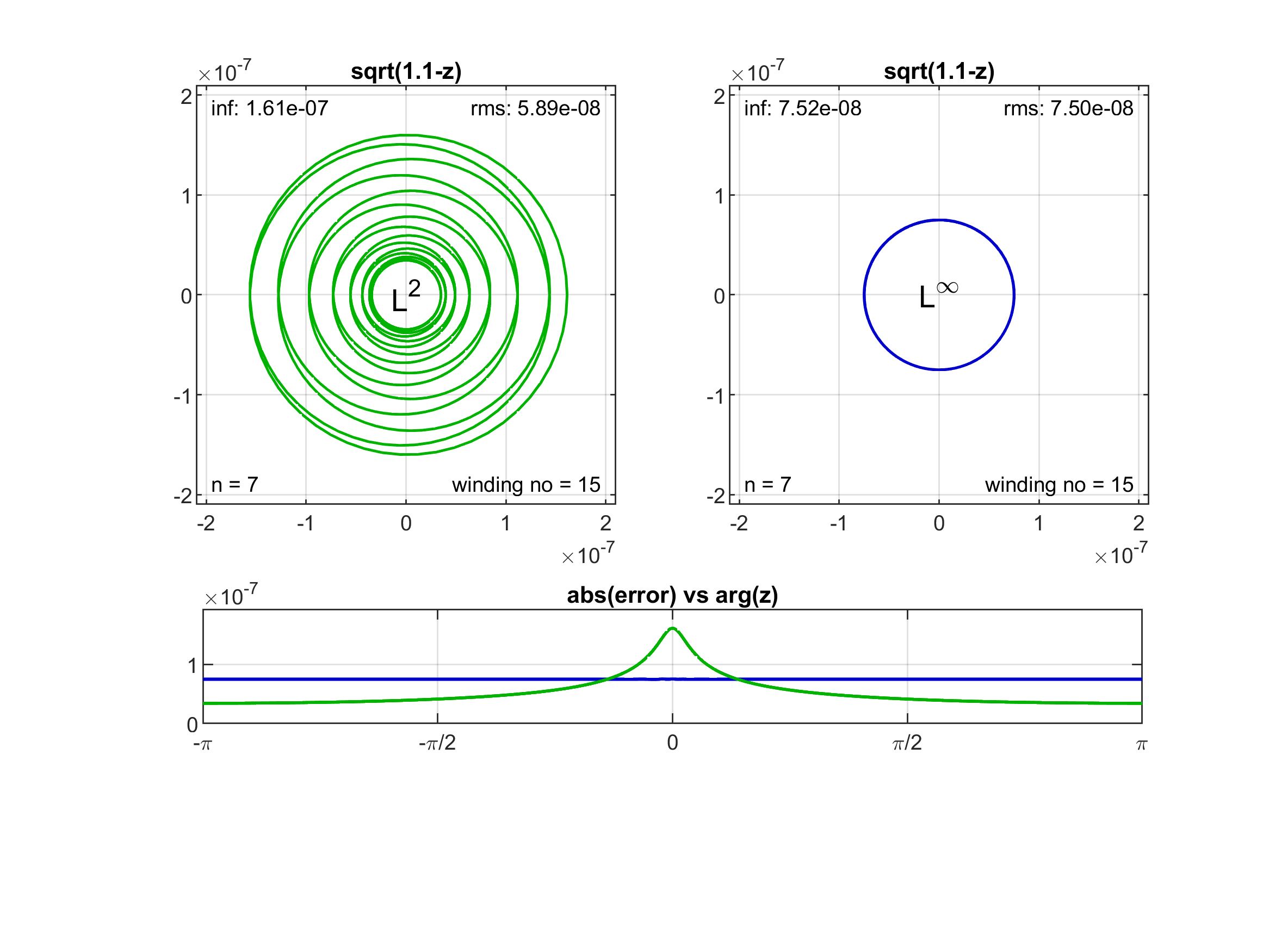}
\end{center}
\caption{\label{fig3}Degree\/ $7$ approximations of\/ $\sqrt{1.1-z}$.}
\end{figure}

In Figure~\ref{fig1} as in many examples, the $L^2$ and $L^\infty$
error curves are topologically the same but quantitatively quite
different.  Both curves wind around the origin $2\kern .3pt n+1$
times, which is the generic number for degree $n$ approximation
though not the universal one.  The $L^2$ curve is roughly circular,
whereas the $L^\infty$ curve is so close to circular that one cannot
see in the plot that it is not exactly so.  The near-circularity
phenomenon was investigated in \cite{nearcirc}, and in this example
the modulus $|e(z)|$ deviates from a constant value by about 3 parts
in $1000$.  Theorem 6.3 of \cite{nearcirc} indicates that for $f(z)
= \exp(z)$ instead of $\exp(4z)$, the deviation from circularity
would be less than one part in a million.

In the upper corners, Figure~\ref{fig1} lists the errors $\|e\|_2^{}$ and
$\|e\|_\infty^{}$ for both the $L^2$ and $L^\infty$ best approximations.
Our definitions imply that these numbers are ordered as follows:
\begin{displaymath}
\|\et\|_2^{} \le \|\ei\|_2^{} \le \|\ei\|_\infty^{} \le \|\et\|_\infty^{}.
\end{displaymath}

As our next examples, Figures 2 and 3 show approximations of $f(z)
= \sqrt{1.1-z}$ of degrees $n=1$ and $7$.  It is interesting to
note here the difference between the plots of error curves, in
the upper part of each figure, and of $|e(z)|$ as a function of
$\hbox{arg}(z)$, in the lower part.  From the error curves, one
might be puzzled as to how the green approximations can be better
in the 2-norm than the blue ones.  The lower plots, however, make
it clear that the $L^2$-optimal approximations have large errors
only on a small portion of $S$ near the singularity at $z=1.1$.

It is also interesting to note in Figures 2 and 3 that, whereas
$L^\infty$ error curves tend to approach perfect circles as
$n\to\infty$, this need not be the case in $L^2$.  Thus in a relative
sense, $\et$ and $\ei$ need not approach each other as $n\to\infty$.
In an absolute sense, however, they do approach each other, for it is
known that $\lim_{n\to\infty}^{} E_2^{1/n}$ and $\lim_{n\to\infty}
E_\infty^{1/n}$ are equal for a wide class of functions $f$
\cite[Thm.~5]{bsy}.  Such results are proved by methods related to
the Hermite integral and potential functions discussed in section~4.

\begin{figure}
\begin{center}
\includegraphics[trim=30 65 20 10, clip, scale=.93]{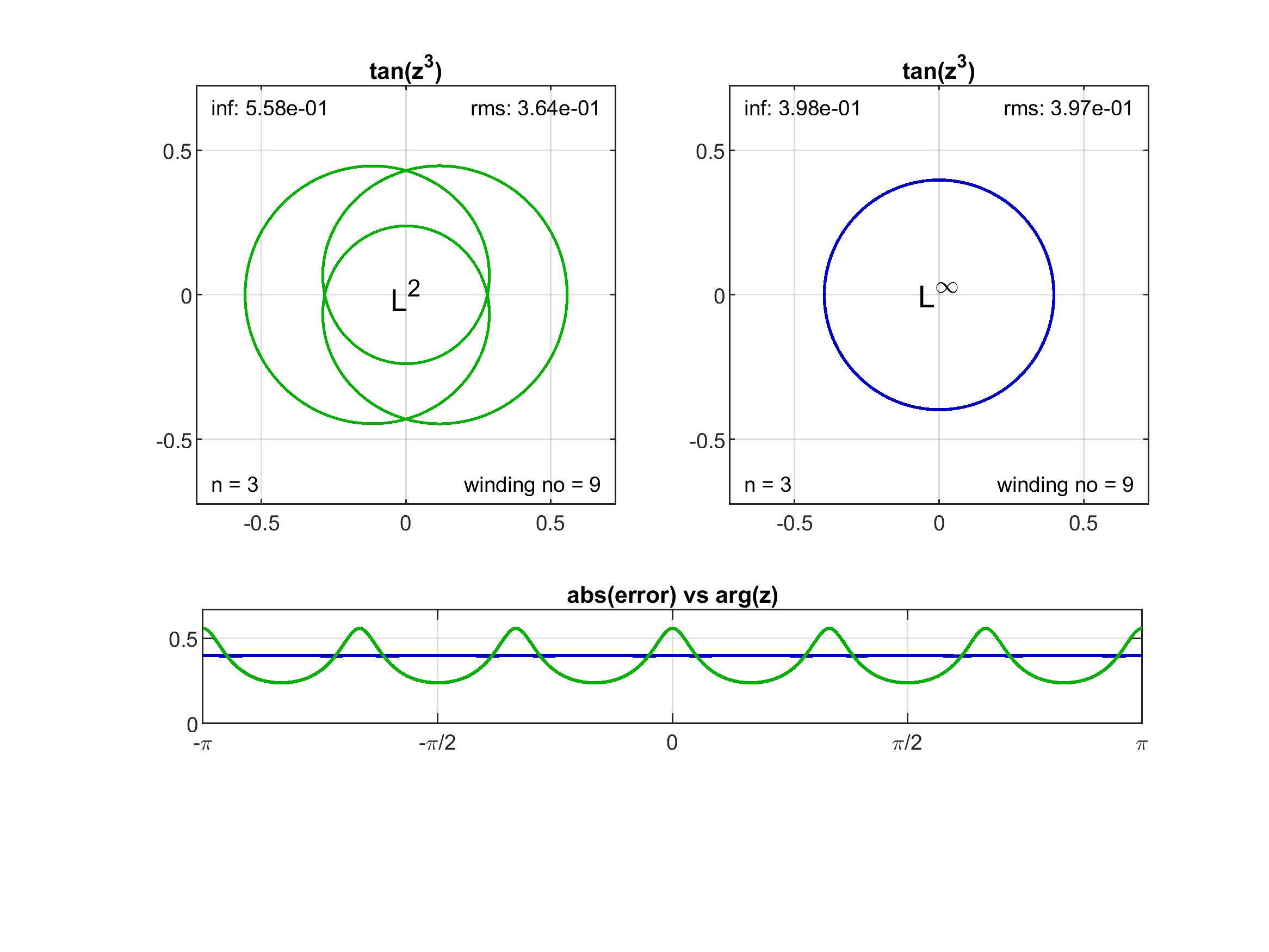}
\end{center}
\caption{\label{fig4}Degree\/ $3$ approximations of\/ $\tan(z^3)$.  The
winding numbers exceed the generic value\/ $2\kern .3pt n+1 = 7$.
The $L^2$ best approximation is nonunique, since it could be rotated by $\pi/3$, and
the unique $L^\infty$ best approximation is a polynomial, $\ri(z) \approx 1.1596z^3$.}
\end{figure}

As a fourth example, Figure 4 shows degree $n=3$ approximations
of $f(z) = \tan(z^3)$, a function with six-fold angular symmetry
about $z=0$.  A new feature here is that the winding number for
both approximations is $\omega = 9$, which exceeds $2\kern .3pt
n+1 = 7$.  The unique $L^\infty$ best approximation is actually
a polynomial, $\ri(z) \approx 1.159500940306 z^3$, so its error
curve is six-fold symmetric.  The best $L^2$ approximation, on the
other hand, has three finite poles equally spaced along a circle,
so it cannot be exactly six-fold symmetric despite the impression
given by the green curve of $|e(z)|$ against $\hbox{arg}(z)$ in the
lower image.  The green error curve in the upper-left image appears
to have winding number only 3, but because of the symmetry, it is
actually traversed three times.  This best approximation must be
nonunique, since an equivalent approximation could be obtained by
a rotation by an angle $\pi/3$ with a negation of sign.  We comment
on such degeneracies in the next section.

\section{Optimality conditions}
We make some comments on optimality conditions, without attempting
a full discussion.

For $L^\infty$, as mentioned in the introduction, best approximations
on the unit disk are known to exist \cite{walsh31} but need not
be unique \cite{gt}.\footnote{For example, Theorem 1 of \cite{gt}
implies that $f(z)=z^2+z^5$ has nonunique best approximations of
degree $n=1$.} The best known steps toward characterizing them
are the {\em local Kolmogorov condition\/}, which is necessary for
optimality (local as well as global), and the {\em global Kolmogorov
condition\/}, also known as the {\em Meinardus-Schwedt condition\/},
which is sufficient for optimality \cite{gutknecht,meinardus,thiran}.
We are not aware of a condition that is both necessary and
sufficient, and it is notable that the AAA and AAA-Lawson algorithms,
unlike previous less successful methods for $L^\infty$ approximation
and unlike our algorithms here for $L^2$ approximation, are not
based on an attempt to enforce optimality conditions.

As a practical matter, best and near-best $L^\infty$ approximations
can usually be recognized by their nearly circular error curves, as
illustrated in all the figures so far.  In the simplest case, suppose
$r\in\Rn$ has an error curve $e(S)$ that happens to be a perfect
circle with winding number $\omega \ge 2\kern .3pt n+1$.  Then if
$\rtil$ were a better approximation, the difference $r-\rtil$ would
also have to have winding number $\omega \ge 2\kern .3pt n+1$, which
is impossible since $r-\rtil\in R_{2\kern .3pt n}^{}$: so $r$ must
be a best approximation.  (This argument amounts to an application
of Rouch\'e's theorem to the error functions $f-r$ and $f-\rtil$.)
This specialized conclusion generalizes in two important ways.
First, if $r\in\Rn$ is of degree less than $n$, belonging to $\Rnu$
for some $\nu<n$, then winding number $\omega \ge n+\nu+1$ is enough
to imply optimality.  Second, if the error curve is not perfectly
circular but nearly circular, then the same reasoning ensures that
$r$ is a {\em near-best\/} approximant \cite[Prop.~2.2]{nearcirc}.
This is a generalization to the unit disk of a familiar lemma in
real approximation on an interval associated with the name of de
la Vall\'ee Poussin~\cite{meinardus}:

\medskip
{\em{\sc Theorem 1.  
Nearly circular $\Rightarrow$ Near-best for $L^\infty$ approximations.}
Given a function\/ $f$ analytic in the closed unit disk,
suppose\/ $r\in\Rnu$ for some\/ $\nu\le n$ has an error curve\/ 
$e(S)$ that does not pass through the origin and has winding number\/
$\omega \ge \nu + n + 1$.  Then
\begin{displaymath}
\min_{z\in S} |e(z)| \le \Ei \le \max_{z\in S} |e(z)|.
\end{displaymath}
}
\medskip

\noindent Theorem 4.1 of \cite{gutknecht} generalizes this winding
number condition to a related condition of alternation of signs
along the unit circle, which can be used to prove the optimality
of the best $L^\infty$ approximation of Figure~\ref{fig4}.

Note that although error curves for $L^\infty$ best approximations
often look circular to the eye, they can never be exactly circular
unless $f$ is rational.  The proof is that if the error curve is
circular, then by the Schwarz reflection principle, $e = r-f$ can
be analytically continued to the extended complex plane $\complex
\cup \{\infty\}$, making it a meromorphic function with finitely
many poles, hence rational.\footnote{As an undergraduate at Harvard,
the third author spent weeks trying to prove this without success.
Finally he knocked at the office of Prof.\ Ahlfors, who delivered
the argument instantly while standing in the doorway.}

Turning to $L^2$, best approximations again exist but need not
be unique.  A helpful pair of sources on this topic are \cite{bar}
and \cite{bsy}.  In this case there is a necessary condition for best
approximation that is the basis of our computations.  First of all,
a best approximation $\rt\approx f$ must always satisfy $\rt(0) =
f(0)$.  This follows from the property that the $L^2$ norm over $S$
of an analytic function with Taylor series $a_0 + a_1z + \cdots$
is the 2-norm of its coefficient vector $(a_0, a_1, \dots)^T$,
as can be shown by orthogonality; thus if an approximation $r$
does not interpolate $f$ at $z=0$, it can be improved by addition
of a constant.  In addition, there must be at least $2\kern .3pt n$
further points of interpolation lying in pairs---that is, points
of {\em Hermite interpolation\/}---located at the reflections
$\{\bar\pi_k^{-1}\}$ in the unit circle of the poles $\{\pk\}$
of $\rt$.  The following theorem appears to have been published
first in \cite{Levin69}, generalizing an earlier result restricted
to the case of simple poles \cite{erohin} (in Russian).  Analogous
results for $L^2$ approximation in a half-plane can be found in
\cite{MeierL67}.

\medskip
{\em {\sc Theorem 2.
Interpolation conditions for best $L^2$ approximation.}
If\/ $f\not\in \Rnm$ is analytic in the closed unit disk,
suppose\/ $\rt\in\Rn$ is a best\/ $L^2$ approximation of\/ $f$.
Then\/ $\rt\not\in R_{n-1}^{}$,
so\/ $\rt$ has exactly\/ $n$ finite or infinite poles\/ 
$\pi_1^{},\dots,\pi_M^{}$, counted with multiplicities
$m_1^{},\dots,m_M^{}$; thus $\sum_{k=1}^M \mk = n$.
Moreover,\/ $\rt$ interpolates\/ $f$ at the origin,
\begin{equation}
\rt(0) = f(0),
\label{char0}
\end{equation}
and in at least\/ $2\kern .3pt n$ additional points counted
with multiplicity in the following sense:
\begin{equation}
r_2^{*(\kern .3pt j)}(\bar\pi_k^{-1}) = f^{(\kern .3pt j)}(\bar\pi_k^{-1}), \quad
j = 0,1,\ldots, 2\kern .3pt \mk-1, \quad k=1,2,\ldots,M.
\label{char}
\end{equation}
If one of the poles is infinite, say\/ $\pi_1^{}=\infty$, then
the $k=1$ case of\/ $(\ref{char})$ becomes
\begin{equation}
r_2^{*(\kern .3pt j)}(0) = f^{(\kern .3pt j)}(0), \quad j = 1,2,\ldots, 2\kern .3pt \mk.
\label{charinf}
\end{equation}
}
\medskip

\begin{figure}
\begin{center}
\includegraphics[trim=30 65 20 10, clip, scale=.93]{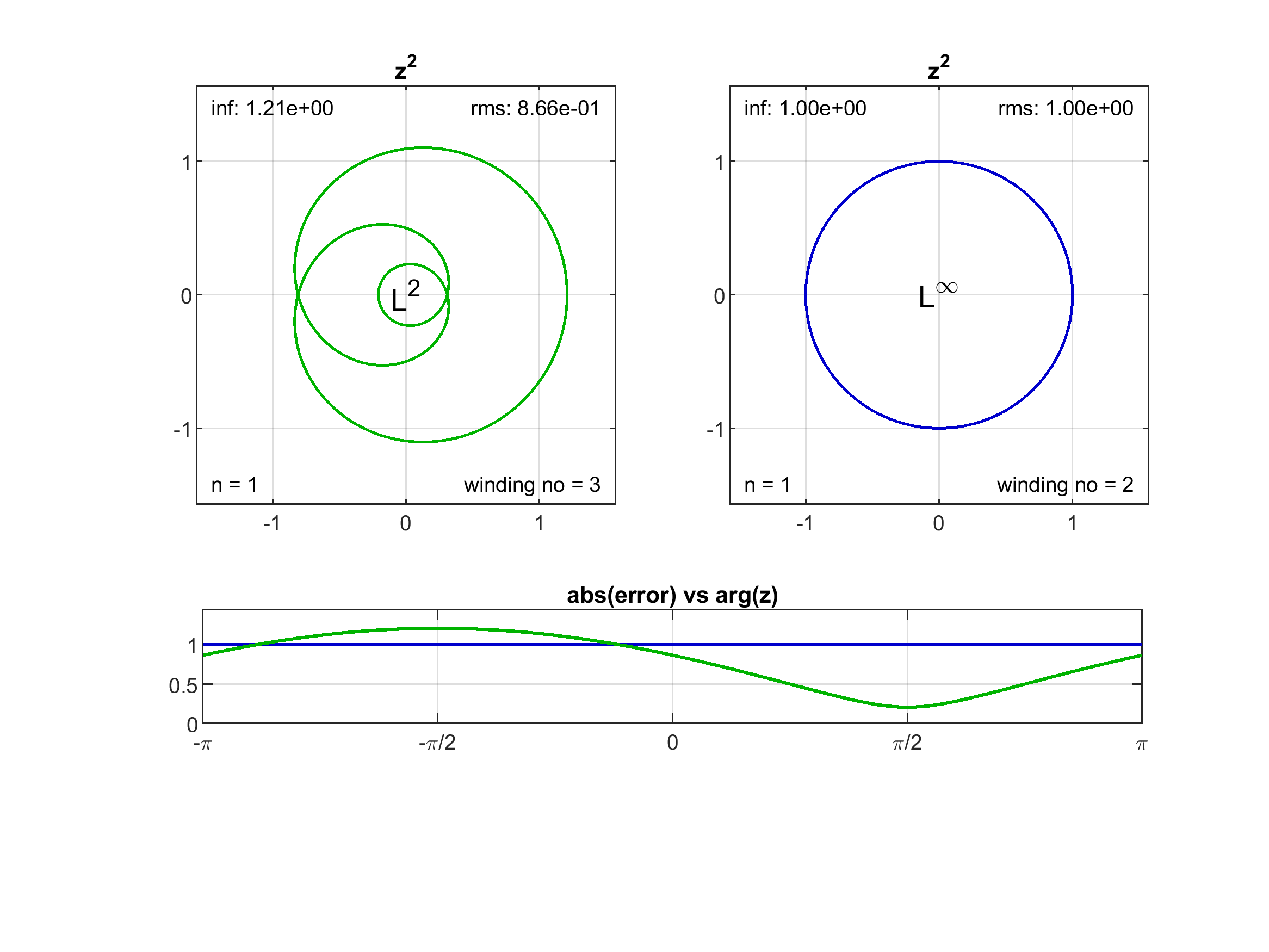}
\end{center}
\caption{\label{fig5}Degree\/ $1$ approximations of\/ $z^2$.  In $L^2$,
this function has a continuum of nonunique best approximations,
as described in the text, which have more points
of interpolation than the (unique) $L^\infty$ best approximant.  The $L^\infty$ error
curve is our first example that is exactly rather than just
nearly circular.}
\end{figure}

Theorems 3.1 and 3.2 imply an interesting difference between
$L^\infty$ and $L^2$ approximation.  For $L^\infty$, the best
approximation in $\Rn$ may be ``defective,'' belonging to a space
$\Rnu$ with $\nu< n$, but for $L^2$ it must always have exact
degree $n$, implying that $\Et$ is a strictly decreasing sequence
as $n\to\infty$ whenever $f$ is not itself rational.  As the
simplest possible example, illustrated in Figure~\ref{fig5},
consider best degree $n=1$ approximation of $f(z) = z^2$
\cite[p.~272]{Levin69}.  In $L^\infty$, the best approximation is
$\ri(z) = 0$, with $\Ei=1$ and winding number $2$, but in $L^2$,
$\rt(z) = (ze^{i\theta}/2)/(ze^{-i\theta}-\sqrt2 \kern .3pt i)$ is
a best approximation for any $\theta$, with $\Et = \sqrt 3/2 \approx
0.866$ and winding number $3$.  Thus the best $L^2$ approximation is
nonunique, as must hold more generally in any case where $f$ is even
and $n$ is odd \cite{bar}.  If $f$ is perturbed slightly in this
example, then instead of a nonunique global best approximation,
we can expect to find that there may be a unique global best
approximation and also one or more distinct approximations that
are locally but not globally optimal~\cite{bar86}.  As mentioned in
the last section, the $L^2$ best approximation of Figure~\ref{fig4}
is also nonunique.

\section{Hermite integral and potential functions}

The theory of rational approximations of an analytic function $f$
on a complex domain $K$ has been developed extensively in the
past century.  The fundamental tool is a Hermite contour integral
involving a potential function $\phi(z)$.  For a sharp form of
the estimate applicable directly to the generic situation, we
assume there are at least $2\kern .3pt n+1$ interpolation points
$\zeta_0,\dots,\zeta_{2\kern .3pt n}\in K$ with $r(\zk) = f(\zk)$,
and we also assume that $r$ has $n$ poles, some of which may
be at $\infty$.  The poles and interpolation points are counted
with multiplicity, and in the case of $L^2$ best approximation,
Theorem~2 asserts that we can expect many interpolation
points of multiplicity $2$.\ \ As a physical interpretation of the
potential function (more precisely, the log of its absolute value),
we imagine that each interpolation point corresponds to a positive
point charge of strength $1$, and each pole is a negative point
charge of strength $2$.  Here is the theorem as stated in section
12 of \cite{acta}, a special case of \cite[Lemma 2]{stahl89},
with roots in earlier work of Walsh and Gonchar.

\medskip
{\em {\sc Theorem 3.}
Let\/ $f$ be analytic in the closure of a Jordan region\/ $\Omega$
bounded by a Jordan contour\/ $\Gamma$, and let the
degree\/ $n$ rational function\/ $r(z)$
interpolate\/ $f$ in\/ $2\kern .3pt n+1$ points
$\zeta_0,\dots,\zeta_{2\kern .3pt n}$ of a compact set $K\subset \Omega$,
counted with multiplicity.
Let\/ $r\in\Rn$ have $n$ finite poles $\pi_1^{},\dots,\pi_n^{}$ outside $\Omega$, and define
\begin{equation}
\label{potfun}
\phi(z) = \prod_{k=0}^{2\kern .3pt n} (z-\zeta_k)
\! \left/ \prod_{k=1}^n (z-\pi_k)^2\right..
\end{equation}
Then for any $z\in K$, 
\begin{equation}
\label{hermite}
f(z) - r(z) = {1\over 2\pi i}
\int_\Gamma {\phi(z)\over \phi(t)}\kern 1pt {f(t) \over t-z} \kern 1pt dt.
\end{equation}
If\/ $r$ has poles at $\infty$, then $(\ref{hermite})$ still
holds with these poles dropped from the quotient $(\ref{potfun})$.
}
\medskip

If {\em all\/} the poles are at $\infty$, then $r$ reduces to
a polynomial and we recover the Hermite integral for polynomial
approximation presented for example in~\cite[Thm.~11.1]{atap}.

The power of Theorem 3---or perhaps we should say the
power of rational approxi\-mation---comes from the behavior of the
quotient $\phi(z)/\phi(t)$ in (\ref{hermite}).  When the negative
charges (poles) are well separated from the positive ones (interpolation
points), this ratio will
be exponentially small on $K$, establishing exponential decrease of
the errors $\Et$ and $\Ei$ as $n\to\infty$.  The theorem applies
to any rational approximation $r\approx f$ with $n$ finite poles
that interpolates $f$ in at least $2\kern .3pt n+1$ points of $K$.
Having this many interpolation points is the generic situation
for best and near-best approximations, and in cases when there
are fewer, one may be able to apply the more general result of
\cite[Lemma 2]{stahl89}, in which only linearized interpolation
after multiplication by a denominator polynomial is required.

\begin{figure}
\begin{center}
\includegraphics[trim=15 135 5  5, clip, scale=.80]{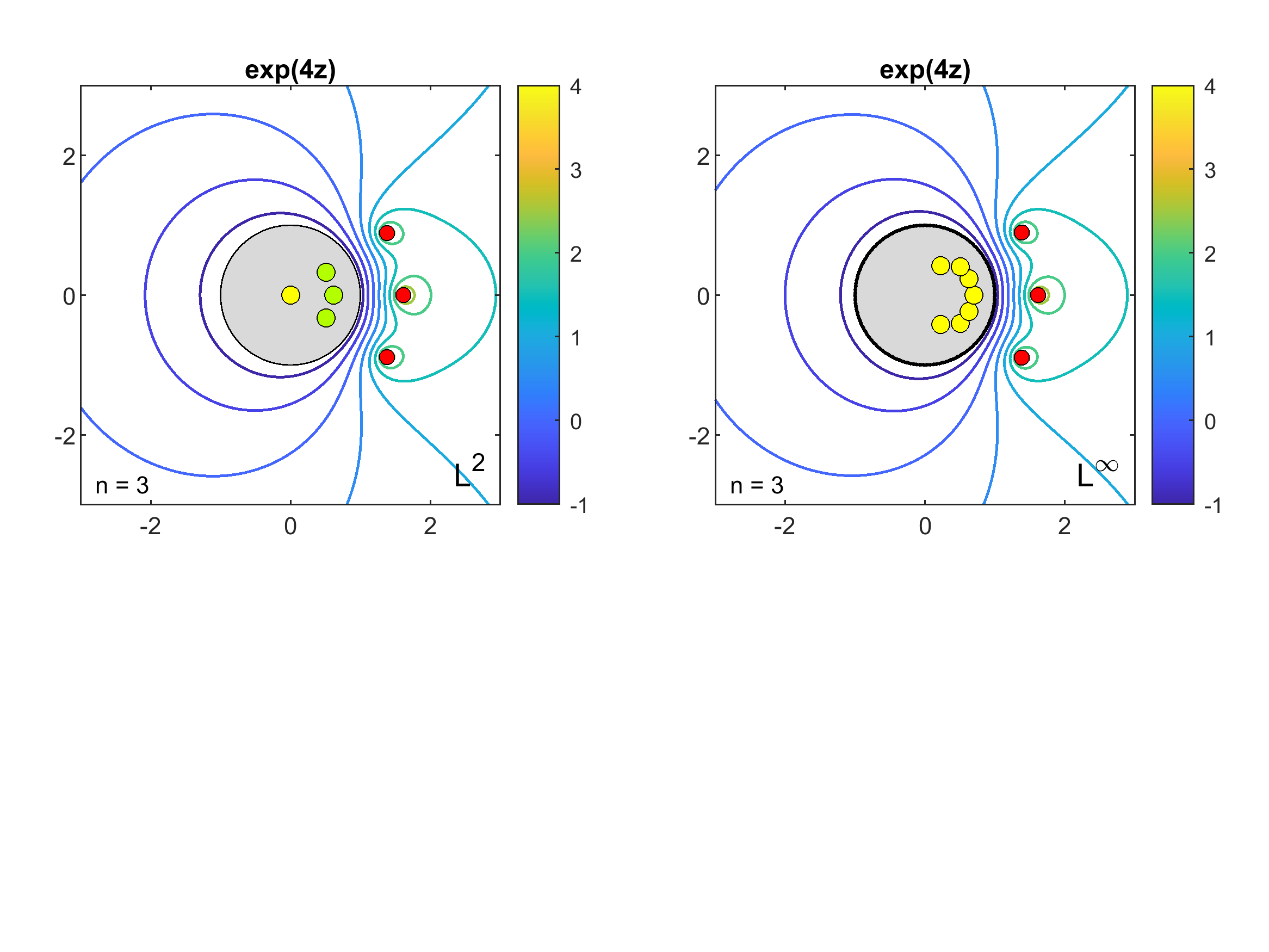}
\end{center}
\vskip -3pt
\caption{\label{fig1p}Potential contours for degree\/ $3$
approximations of\/ $\exp(4z)$, following Figure~$\ref{fig1}$.
Contours show integer and half-integer levels of\/ $\log_{10}|\phi(z)|$, with
$\phi$ defined by\/ $(\ref{potfun})$.
Red dots mark
poles of the approximant $r$, yellow dots mark interpolation points where
$r=f$, and green dots marks double (Hermite) interpolation points with
$r=f$ and $r'=f'$.  By Theorem\/~$2$, for $L^2$ approximation
there is a yellow dot at the origin and
green dots at the
points of reflection of the red dots in the unit circle.}
\end{figure}

\begin{figure}
\begin{center}
\includegraphics[trim=15 135 5  5, clip, scale=.80]{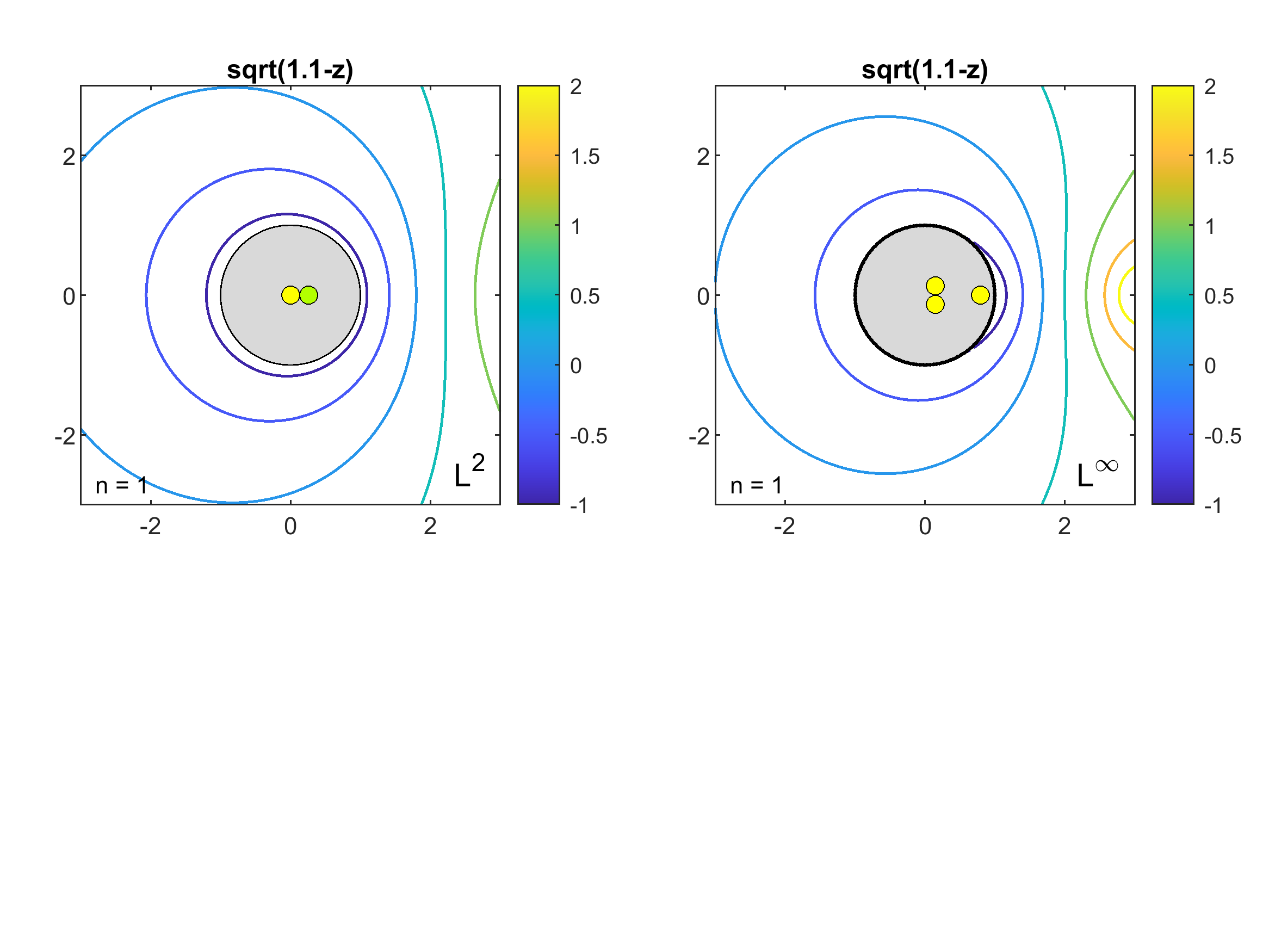}
\end{center}
\vskip -3pt
\caption{\label{fig2p}Potential contours for degree\/ $1$ approximations
of\/ $\sqrt{1.1-z}$, following Figure~$\ref{fig2}$, with dots
as in Figure\/~$\ref{fig1p}$.  The pole $\pz$ is off-scale to the right at
$\pz\approx 3.865$ for $L^2$ and\/ $\pz\approx 3.146$ for $L^\infty$.}
\end{figure}

\begin{figure}
\begin{center}
\includegraphics[trim=15 135 5  5, clip, scale=.80]{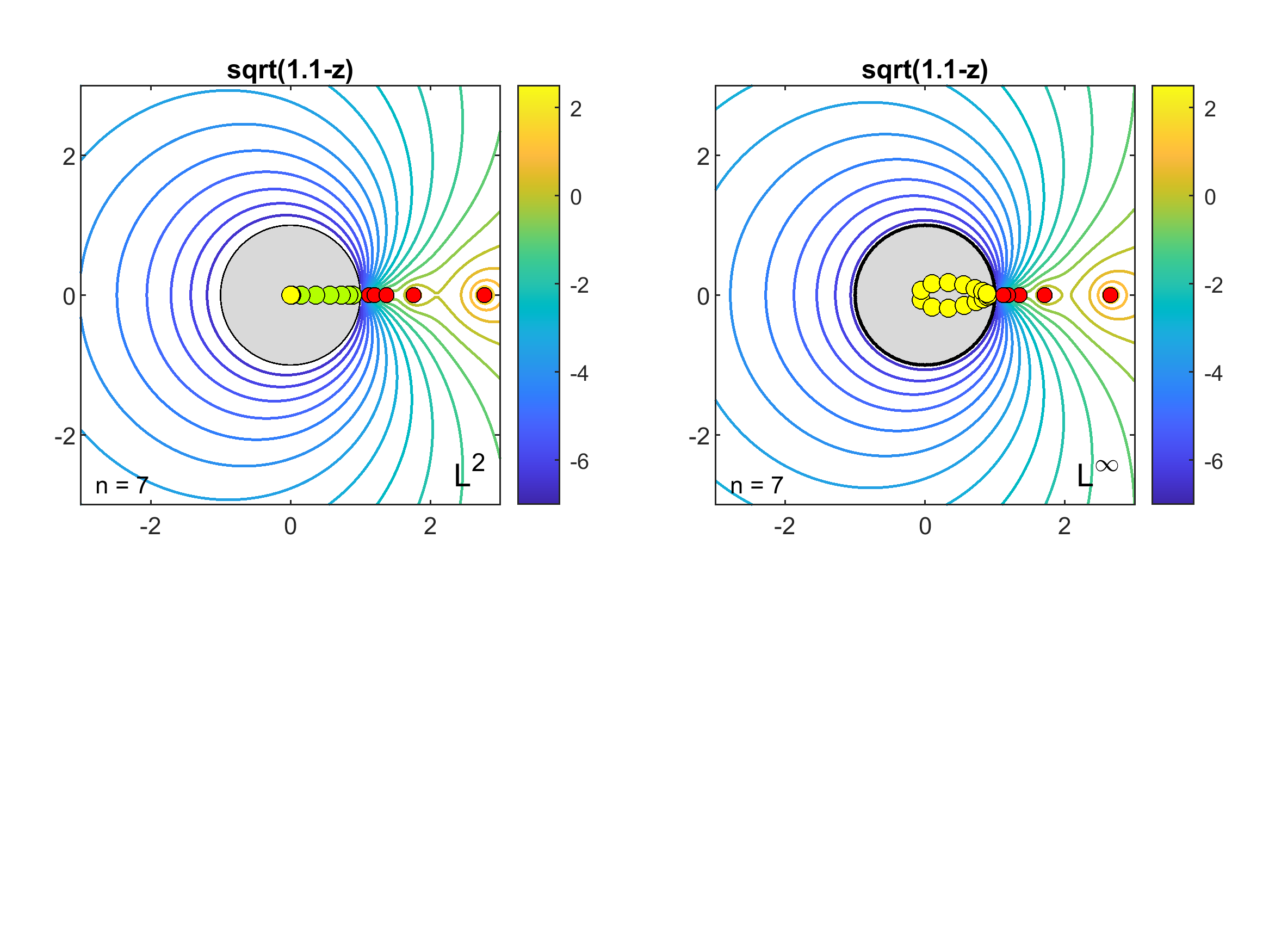}
\end{center}
\vskip -3pt
\caption{\label{fig3p}Potential contours for degree\/ $7$ approximations
of\/ $\sqrt{1.1-z}$, following Figure~$\ref{fig3}$.  In each case two of
the seven poles $\{\pk\}$ are off-scale to the right.}
\end{figure}

\begin{figure}
\begin{center}
\includegraphics[trim=15 135 5  5, clip, scale=.80]{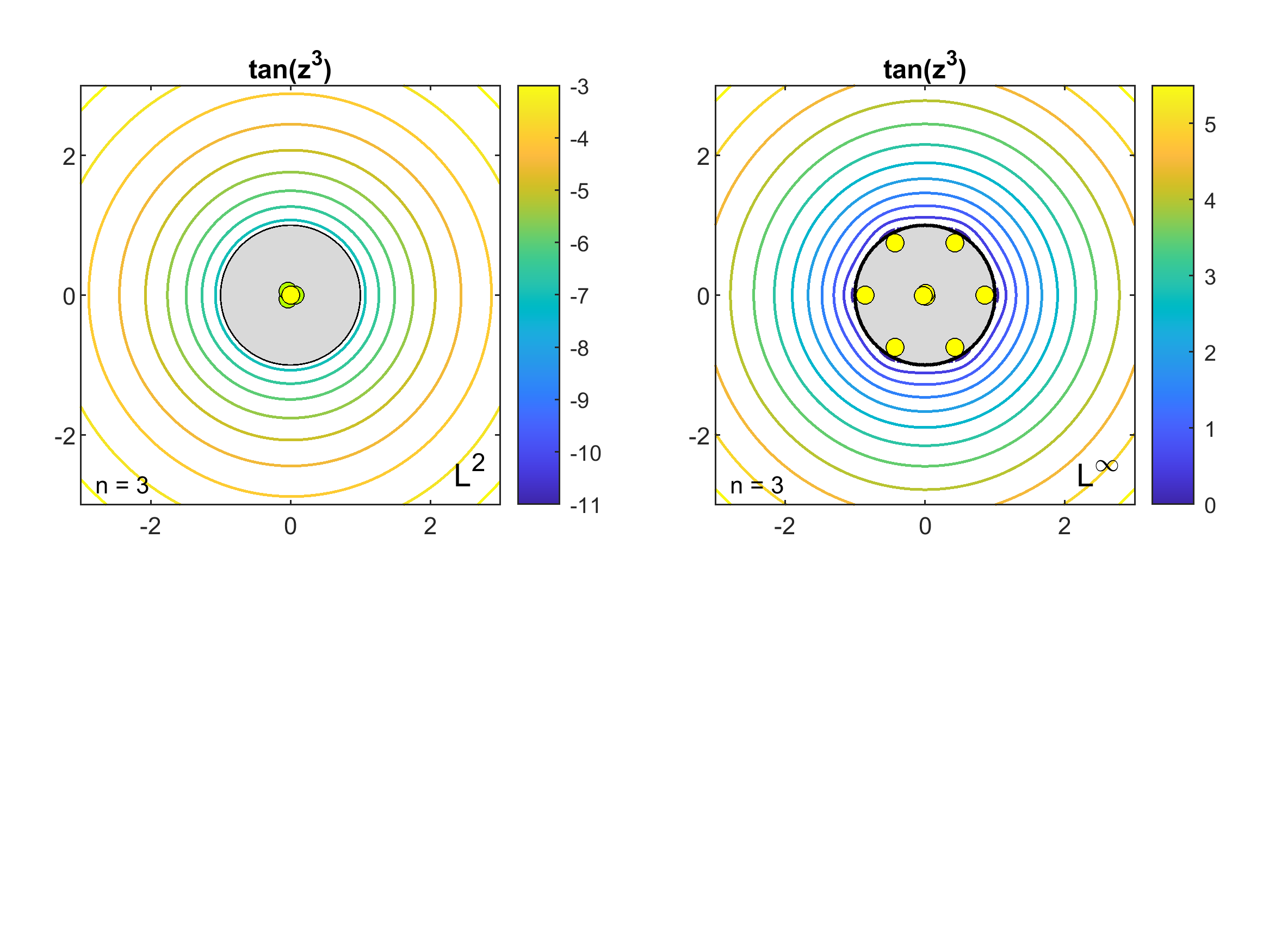}
\end{center}
\vskip -3pt
\caption{\label{fig4p}Potential contours for degree\/ $3$ approximations
of\/ $\tan(z^3)$, following Figure~$\ref{fig4}$.  The poles for\/
$L^2$ approximation are off-scale
with modulus about\/ $15.95$, leading to three Hermite interpolation points
with modulus about\/ $0.06$ and three additional interpolation points at the
origin.  For $L^\infty$ the approximation is
a polynomial with a triple pole at $\infty$.  Note that neither of these
images shows $6$ strings of poles radiating out from the disk to match the poles
of\/ $\tan(z^3)$ as one might expect, for the degree $n=3$ is too small.}
\end{figure}

Figures \ref{fig1p}--\ref{fig4p} show contour plots of $|\phi(z)|$
for the examples of Figures \ref{fig1}--\ref{fig4}.  In each plot,
red dots mark poles of $r(z)$, yellow dots mark interpolation points
with $r(z) = f(z)$, and green dots mark Hermite interpolation
points associated with Theorem~2 with $r(z) =
f(z)$ and $r'(z) = f'(z)$.  The contours are half-integer
level curves of $\log_{10}(|\phi(z)|)$, and the first thing
to note in the figures is that in each case the unit circle,
the boundary of our approximation domain $K$, is approximately a
level curve.  This correlates with near-optimality, since the error
(\ref{hermite}) will be small when the points $\zk\in K$ and $\pk\in
\complex\backslash\overline{\Omega}$ are in a near minimal-energy
configuration.  In the case of best $L^2$ approximation, by
Theorem~2, the unit circle is {\em exactly\/} a level
curve, for $\phi$ is a multiple of a finite Blaschke product (a
rational function whose poles and zeros are reflections of each
other in the unit circle).

Looking first at Figure~\ref{fig1p}, we see a fine illustration of
the generic situation for $L^2$ and $L^\infty$ best approximations.
Both the $L^2$ and $L^\infty$ approximants have three poles to the
right of the unit circle, which lie at nearly but not exactly the
same locations in the two cases.  This approximate agreement of
$L^2$ and $L^\infty$ poles is common in our experience.  The $L^2$
best approximation shows the poles reflected as three green dots,
with double interpolation at each; in addition the yellow dot at
the origin marks the seventh interpolation point.  The $L^\infty$
best approximation also has seven interpolation points in the disk,
but now they are distinct and nonzero.

Figure~\ref{fig2p} is similar with $n=1$, but perhaps less clear than
Figure~\ref{fig1p} since the poles are off-scale.  Figure~\ref{fig3p}
shows the configuration with the degree increased to $7$, with the
poles clustered near the singularity at $z=1.1$.

\section{\boldmath$L^2$ approximation algorithm: state space form}

The starting point of our $L^2$ computations is the
Transfer Function (TF)
variant of the Iterative Rational Krylov Algorithm (IRKA)~\cite{IRKA},
called TF-IRKA~\cite{bg}.
IRKA was originally developed as a method for the
$H^2$-optimal model order reduction of large-scale linear dynamical systems
in continuous time.  Thus, the approximation domain
was the right half-plane;
a closely related algorithm set on
the disk exterior is called MIRIAm~\cite{bunse}.
The input-to-output response of such a
system with $N$ ordinary differential equations,
$p$ input ports, and $q$ output ports is characterized
by a degree-$N$ matrix-valued rational function of
dimensions $q\times p$, the \emph{transfer function,}
\begin{equation}
    \label{eq:TransferFunction}
    \H(z) = \C^{\top}(z\kern .5pt \E-\A)^{-1}\B + \D,
\end{equation}
where $\E\in \RNN$ (nonsingular), $\A \in \RNN$, $\B \in \RNp$, $\C
\in \RNq,$ and $\D \in \Rqp$ are the \emph{state-space matrices} of
$\H(z)$; see~\cite{Antbook} for a general resource on system-theoretic
model reduction.  We refer to~\eqref{eq:TransferFunction} as the
state-space form of the rational function, and any degree $N$
matrix-valued rational function of dimensions $q\times p$ can be
expressed in this form~\cite[sec.~2.6.1]{dullerud}.  (TF-)IRKA is an algorithm for computing low
degree rational approximants of $\H$ that satisfy a set of Hermite
interpolatory optimality conditions, akin to~\eqref{char}\footnote{In
its original setting of~\cite{bg,IRKA}, the domain of analyticity
is the open right half-plane, and the optimal interpolation
points are the reflections of the poles in the imaginary axis,
i.e., $\{-\overline{\pi}_k\}$.}.  In the context of model order
reduction, the computed $\rt$ is the transfer function of another
linear dynamical system of degree $n$ with the same structure and
dimensions $q\times p$ as~\eqref{eq:TransferFunction}.

For the following discussion, assume that the poles are simple,
distinct, and finite.  (TF-)IRKA computes
$\rt$ by an iterative process that adjusts
the poles at each step until a fixed point is reached.
Specifically, at
each step, it produces a new rational approximant by interpolating doubly at
the reflected poles of the previous iterate.  Upon convergence---once
the poles of the rational approximant stop changing up to a
tolerance---the optimality conditions~\eqref{char} are satisfied.
Note that although these conditions are only necessary for
optimality, since they will hold at maxima and stationary points as
well as minima, in practice in our iteration they can be expected to
behave as necessary and sufficient conditions for a local minimum.
The reason is that stationary points other than minima act as
repellers of the IRKA iteration, so the algorithm is unlikely to
converge to them.  (This repelling property has been published as
a theorem in the case of $L^2$ approximation on a half-plane in
the appendix of \cite{kraj} and presumably holds on the disk too.)

The difference between
IRKA and its TF variant is in how it achieves the interpolation.
In its traditional formulation,
at each step, IRKA computes a rational interpolant by projecting
the matrices $\E, \A, \B, \C$ onto particularly chosen rational
Krylov subspaces.
By contrast, TF-IRKA is data-driven, meaning that it only requires
the ability to sample $f$ and $f'$ (not necessarily rational) at points in
$\complex$, and these data are used to construct $\rt$ in the
form~\eqref{eq:TransferFunction}.  This is accomplished using the
\emph{Loewner-matrix} approach for interpolation.
For the remainder of this section, we switch from
TF-IRKA's original formulation to describe how we utilize
it to compute best $L^2$ approximations to functions analytic in the
disk.  That is, the final approximation will satisfy both 
(\ref{char0}) and (\ref{char}).
For our scalar approximation problems, we take $p=q=1$.
Given samples $f_i=f(\sigma_i)$ and derivatives
$f_i'=f'(\sigma_i)$ at the points $\sigma_i\in\complex$,
define the Hermite Loewner and shifted Hermite Loewner
matrices $\L\in\complex^{n\times n}$ and $\M\in\complex^{n\times n}$ by
\begin{equation}
\label{eq:LoewnerMatrices}
\L_{ij}=\begin{cases}
  \displaystyle{\frac{f_i-f_j}{\sigma_i-\sigma_j}} &\mbox{if}~~i\neq j\\[8pt]
  f_i' & \mbox{if}~~i=j
\end{cases},~~~\M_{ij}=\begin{cases}
  \displaystyle{\frac{\sigma_if_i-\sigma_jf_j}{\sigma_i-\sigma_j}} &\mbox{if}~~i\neq j\\[8pt]
    f_i+\sigma_if_i' & \mbox{if}~~i=j
\end{cases}
\end{equation}
as well as the data matrix $\Y=\begin{bmatrix}
        f_1&\cdots&f_n
    \end{bmatrix}^{\top}.$
Construct $\E$, $\A$, $\B$, and $\C$ as
\begin{displaymath}
    \E=-\L,~~~\A=-\M,~~~\B=\Y,~~~\C=\Y.
\end{displaymath}
This achieves the $2\kern .3pt n$ interpolation conditions
\eqref{char}. For the other condition (\ref{char0}), define
\begin{displaymath}
\D={d_0\over d_1d_2 + d_3d_0},
\end{displaymath}
where
\begin{displaymath}
d_0=f(0)+\C^{\top}\kern -1.3pt \A^{-1}\kern 1pt\B,\quad
d_1=-\one^{\top}\kern -1.3pt \A^{-1}\kern 1pt \B - 1
\end{displaymath}
and
\begin{displaymath}
d_2=-\C^{\top}\kern -1.3pt \A^{-1}\one - 1, \quad
d_3 = - \one^{\top}\kern -1.3pt \A^{-1}\kern .7pt\one,
\end{displaymath}
and update the $\A$, $\B$, and $\C$ matrices according to
\begin{displaymath}
\A\gets\A + \D, \qquad \B\gets\B-\D, \qquad \C\gets\C-\D,
\end{displaymath}
where $\one\in\R^{n}$ is the vector of all ones.
By construction, the computed interpolant $r(z)$
constructed as in~\eqref{eq:TransferFunction} satisfies
the $2\kern .3pt n+1$ interpolation conditions
\begin{equation*}
f(\sigma_k)=r(\sigma_k),~~f'(\sigma_k)=r'(\sigma_k),~~f(0)=r(0)
\end{equation*}
for $k=1,\ldots,n$.

Due to its origins in model order reduction, (TF-)IRKA has primarily
been employed for computing low-degree rational approximations to
high-degree rational functions.  However, since the iteration is
data-driven, it can equally be employed to compute $L^2$-optimal
approximants to arbitrary functions analytic on a closed half-plane
or disk or disk exterior.  For further details on TF-IRKA and the
Hermite Loewner framework, see~\cite{abgbook}.

\section{\boldmath$L^2$ approximation algorithm: barycentric form}
In this section we provide an alternative
method for computing the rational interpolants at each step of the TF-IRKA iteration.
Recall that the interpolatory barycentric form of a degree $n$ rational function is
\begin{equation}
\label{eq:bcForm}
r(z) = \sum_{i=1}^{n+1}\frac{w_if_i}{z-t_i}\kern -1pt \left/ \kern 2pt
\sum_{i=1}^{n+1}\frac{w_i}{z-t_i}\right.,
\end{equation}
where $\{t_i\}$ are support points, $\{f_i\}$ are function values,
and $\{w_i\}$ are barycentric weights.
For any choice of nonzero weights, $r(t)\to f(t_i)$ as $t\to t_i$, and we
remove the removable singularities by defining
$r(t_i) = f(t_i)$ accordingly.
This interpolatory property of $r$
allows us to satisfy the $n+1$ Lagrange interpolation conditions
of Theorem 2 ``for free'' by taking $t_1^{},\dots,t_n^{}$ as
the reflections of the current poles $\{\pi_i\}$ in the unit
circle and $t_{n+1}^{}=0$.
We then use the weights to satisfy the remaining $n$ Hermite conditions at
$t_1^{},\dots,t_n^{}$.
A value $d_k \in \complex$ is the derivative of a rational function
in barycentric form \eqref{eq:bcForm} at $z=t_k$ if
\begin{displaymath}
    0 = w_kd_k + \sum_{i\neq k}^{n+1} \frac{w_i(f_k-f_i)}{t_k-t_i},
\end{displaymath}
as is readily derived from the Schneider--Werner formula
for the derivative of a bary\-centric quotient~\cite{SchWer}.
Hence, we can enforce that the derivative of the barycentric form
takes prescribed values $\{d_i\}$ at $t_1,\dots, t_n$ by requiring
that the vector of weights $\mathbf w \in \complex^{n+1}$ satisfies
\begin{displaymath}
    \mathbf w \in \textrm{null}(\widehat\L)
\end{displaymath}
where $\widehat\L \in \complex^{n \times (n+1)}$ is the rectangular
Hermite Loewner matrix defined by deleting the last row of $\L$
in \eqref{eq:LoewnerMatrices}.  Similar expressions appear in
\cite{GosInternal,WillMitchell}, but to our knowledge this is
the first time the conditions for an exact Hermite interpolating
rational function in barycentric form have appeared.

For the computations of this paper, we have used the barycentric
method of this section rather than the Loewner framework method
of the last section.  While we have not conducted extensive tests,
we have not noticed significant practical differences between the
two approaches.

\section{Discussion}
Sample Matlab codes used to generate the figures of this paper
can be obtained from the authors.  We have not currently
developed robust software for such computations.

On regions other than disks or half-planes, the AAA-Lawson method
for $L^\infty$ still applies, but the $L^2$ methods of sections 5
and 6 do not.  
However, as mentioned in the introduction, recent work by Borghi
and Breiten represents steps in this direction \cite{borghi24,borghi}.

Rational approximations are not only useful for rational
approximation.  They can play a role in many other computational
problems, and for a survey of dozens of such applications, see
\cite{acta}.  A theme in applications is that rational functions
are often much more powerful than polynomials, especially when
dealing with unbounded or nonconvex domains, or with functions
having singularities on the domain boundary or nearby.  On the
other hand, we are not aware that there is a great difference in
the approximation power of $L^2$ and $L^\infty$ approximations on
bounded domains.

\section*{Data availability statement}
This paper does not make use of any data.

\section*{Acknowledgments}
This project originated in April 2025 in Banff, Alberta at the BIRS
Workshop on Challenges, Opportunities, and New Horizons in Rational
Approximation organized by Anil Damle, Serkan Gugercin, and Heather
Wilber, where MSA and SR met LNT and first began to compare $L^2$
and $L^\infty$.  We are also grateful for helpful suggestions from
Laurent Baratchart and Maxim Yattselev.


\begin{thebibliography}{1}

\bibitem{Antbook}
{\sc A. C. Antoulas\/},
{\em Approximation of Large-Scale Dynamical Systems},
SIAM, 2005.

\bibitem{abgbook}
{\sc A. C. Antoulas, C. Beattie, and S. Gugercin\/},
{\em Interpolatory Methods for Model Reduction},
SIAM, 2020.

\bibitem{GosInternal}
{\sc L. Balicki, I.~V. Gosea, and S. Gugercin,}
{\em Including derivatives to the barycentric representation of rational
interpolants,} unpublished report (2022).

\bibitem{bar86}
{\sc L. Baratchart,}
{\em Existence and generic properties of $L^2$ approximants
for linear systems,} 
IMA J. Math.\ Control \& Infor., 3 (1986), pp.~89--101. 

\bibitem{bar}
{\sc L. Baratchart,}
{\em A remark on uniqueness of best rational approximants
of degree $1$ in $L^2$ of the circle,}
Elect.\ Trans.\ Numer.\ Anal., 25 (2006), pp.~54--66. 

\bibitem{bsy}
{\sc L. Baratchart, H. Stahl, and M. Yattselev,}
{\em Weighted extremal domains and best rational approximation\/},
Adv. Math., 209 (2012), pp.~357--407.

\bibitem{bg}
{\sc C. Beattie and S. Gugercin,}
{\em Realization-independent ${\cal H}_2$ approximation\/},
in 2012 IEEE 51st IEEE Conference on Decision and Control (CDC), IEEE,
pp.~4953--4958.

\bibitem{borghi24}
{\sc A. Borghi and T. Breiten,}
{\em Generalizing the optimal interpolation points for IRKA,}
in MATHMOD Short Contribution volume 2025: 11th Vienna
International Conference on Mathematical Modeling (2025), pp.~7--8.

\bibitem{borghi}
{\sc A. Borghi and T. Breiten,}
{\em Data-driven optimal approximation on Hardy spaces in simply connected domains,}
Adv.\ Comput.\ Math., 51 (2025), pp.~59-1--59-33.

\bibitem{bunse}
{\sc A. Bunse-Gerstner, D. Kubali\'nska, G. Vossen, and D. Wilczek,}
{\em $h_2$-norm optimal model order reduction for large
scale discrete dynamical MIMO systems,}
J. Comput.\ Appl.\ Math., 233 (2010), pp.~1202--1216.

\bibitem{continuum}
{\sc T. A. Driscoll, Y. Nakatsukasa, and L. N. Trefethen,}
{\em AAA rational approximation on a continuum,}
SIAM J Sci.\ Comput., 46 (2024), pp.~A929--A952.

\bibitem{dullerud}
{\sc G. Dullerud and F. Paganini,}
{\em A Course in Robust Control Theory: A Convex Approach,}
Springer Science \& Business Media (2013).

\bibitem{erohin}
{\sc V. D. Erohin,}
{\em On the best approximation of analytic functions by rational
functions with free poles,}
Dokl.\ Akad.\ Nauk SSSR, 128 (1959), pp.~29--32. (Russian)

\bibitem{grimm}
{\sc J. Grimm,}
{\em Rational approximation of transfer functions in the Hyperion software,}
Technical Report 4002, INRIA (2000).

\bibitem{IRKA}
{\sc S. Gugercin, A. C. Antoulas, and C. Beattie\/},
{\em ${\cal H}_2$ model reduction for large-scale linear dynamical systems\/},
SIAM J. Matrix Anal.\ Applics., 30 (2008), pp.~609--638.

\bibitem{gutknecht}
{\sc M. H. Gutknecht,}
{\em On complex rational approximation part I: The characterization
problem,} Computational Aspects of Complex Analysis:
Proceedings of the NATO Advanced Study Institute held at Braunlage,
Harz, Germany, July 26–August 6, 1982, Springer Netherlands (1983).

\bibitem{gt}
{\sc M. H. Gutknecht and L. N. Trefethen,}
{\em Nonuniqueness of best rational approximations on the
unit disk,} J. Approx.\ Th., 39 (1983), pp.~275--288.

\bibitem{kraj}
{\sc W. Krajewski, A. Lepschy, M. Redivo-Zaglia, and U. Viaro,}
{\em A program for solving the $L_2$ reduced-order model problem
with fixed denominator degree,}
Numer.\ Algs., 9 (1995), pp.~355--377.

\bibitem{Levin69}
{\sc A. L. Levin,}
{\em The distribution of poles of rational functions
of best approximation and related questions,}
Math.\ USSR-Sbornik, 9 (1969), pp.~267--274.

\bibitem{mayoant}
{\sc A. J. Mayo and A. C. Antoulas\/},
{\em A framework for the solution of the generalized realization problem},
Lin.\ Alg.\ Applics., 425 (2007), pp.~634--662.

\bibitem{MeierL67}
{\sc L. Meier and D. Luenberger,}
{\em Approximation of linear constant systems,}
IEEE Trans.\ Aut. Control, 2 (1967), pp.~585--588.

\bibitem{meinardus}
{\sc G. Meinardus,}
{\em Approximation of Functions: Theory and Numerical Methods,}
Springer (1967).

\bibitem{WillMitchell}
{\sc W. Mitchell,}
{\em Two intriguing variants of the AAA algorithm for rational approximation,}
arXiv:2508.11169 (2025).

\bibitem{aaa}
{\sc Y. Nakatsukasa, O. S\`ete, and L. N. Trefethen,}
{\em The AAA algorithm for rational approximation\/},
SIAM J. Sci.\ Comput., 40 (2018), pp.~A1494--A1522.

\bibitem{lawson}
{\sc Y. Nakatsukasa and L. N. Trefethen,}
{\em An algorithm for real and complex rational minimax approximation\/},
SIAM J. Sci.\ Comput., 42 (2020), pp.~A3157--A3179.

\bibitem{acta}
{\sc Y. Nakatsukasa and L. N. Trefethen,}
{\em Applications of AAA rational approximation\/},
Acta Numer., 35 (2026), to appear.

\bibitem{SchWer}
{\sc C. Schneider and W. Werner,}
{\em Some new aspects of rational interpolation,}
Math.\ Comput.\ 47 (1986), pp.~285--299.

\bibitem{stahl89}
{\sc H. Stahl,}
{\em On the convergence of generalized Pad\'e approximants,}
Constr.\ Approx., 5 (1989), pp.~221--240.

\bibitem{thiran}
{\sc J.-P. Thiran and M.-P. Istace,}
{\em Optimality and uniqueness conditions in complex
rational Chebyshev approximation with examples,}
Constr.\ Approx., 9 (1993), pp.~83--103.

\bibitem{nearcirc}
{\sc L. N. Trefethen,}
{\em Rational Chebyshev approximation on the unit disk\/},
Numer.\ Math., 37 (1981), pp.~297--320.

\bibitem{atap}
{\sc L. N. Trefethen,}
{\em Approximation Theory and Approximation Practice, extended edition,}
SIAM (2019).

\bibitem{walsh31}
{\sc J. L. Walsh,}
{\em The existence of rational functions of best approximation\/},
Trans.\ AMS, 33 (1931), pp.~668--689.

\end{thebibliography}
\end{document}